\documentclass[11pt]{article}
\usepackage{amsmath}
\usepackage{amssymb}
\usepackage{xcolor}
\title{Geronimus transformations for sequences of $d$-orthogonal polynomials}

\begin{document}
\author{D. Barrios Rolan\'{\i}a, J.C. Garc\'{\i}a-Ardila\\
Dpto. Matem\'atica Aplicada a la Ingenier\'{\i}a Industrial\\
Universidad Polit\'ecnica de Madrid }
\date{}
\maketitle
\newtheorem{lem}{Lemma}
\newtheorem{prop}{Proposition}
\newtheorem{defi}{Definition}
\newtheorem{coro}{Corollary}
\newtheorem{rem}{Remark}
\newtheorem{nota}{Note}
\newtheorem{teo}{Theorem}
\def\adots{\mathinner{\mkern2mu\raise1pt\hbox{.}\mkern2mu\raise4pt\hbox{.}\mkern2mu\raise7pt\hbox{.}\mkern1mu}}
\newcommand{\se}[1]{\medskip \begin{center}{\section{#1}} \end{center}
\medskip }

 \begin{abstract}
In this paper an extension of the concept of Geronimus transformation for sequences of $d$-orthogonal polynomials $\{P_n\}$  is introduced. The transformed sequences $\{P^{(k)}_n\}$, for $ k=1,\ldots, d,$ are analyzed and some relationships between these new sequences of polynomials are given. Also the associated Hessenberg matrices and their transformed are studied.
 \end{abstract}

\section{Introduction}
For a fixed $d\in \mathbb{N}$, let us consider a vector of linear functionals $(u_1,\ldots,u_d)$,  
$$
\begin{array}{rclccccc}
u_i: &\mathcal{P}[x]&\longrightarrow & \mathbb{C}\\
 & q(x) & \longmapsto & \langle u_i,q\rangle\,, && i=1,2,\ldots, d\,,
 \end{array}
$$
defined on the space $\mathcal{P}[x]$ of polynomials with complex coefficients. The notion of $d$-orthogonality for this kind of vectors 
$
(u_1,\ldots,u_d)\in \left( \mathcal{P}[x]^\prime\right)^d
$
was introduced in ~\cite{iseghem-1987} and used since 1987 in several applications of Approximation Theory (see \cite{Kalyagin,2018} for example).

\begin{defi}[J. Van Iseghem~\cite{iseghem-1987}]
 We say that the sequence of  polynomials $\{P_n\},\, n\in \mathbb{N},$ is a $d$-orthogonal sequence ($d$-OPS in the sequel) if $\deg P_n(x)=n$ for all $n\in \mathbb{N}$ and there exists a vector of  functionals $(u_1,\ldots,u_d)$  such that
\begin{equation}
\left.
\begin{array}{rclc}
\label{d-rec}
\langle{u_j,x^mP_n}\rangle& = & 0, \quad n\geq md+j, \quad m\geq 0, \\
\langle{u_j,x^nP_{dn+j-1}}\rangle & \ne  & 0,\quad  n\geq 0,
\end{array}	
\right\}
\end{equation}
 for each $j=1,\ldots, d$.
 
We say that the vector of  functionals $(u_1,u_2,\ldots, u_d)$ is regular if there exists a sequence of polynomials  $\{P_n\},\, n\in \mathbb{N},$ satisfying \eqref{d-rec}.
\end{defi}

\begin{defi}
If \eqref{d-rec} is verified, then we say that $(u_1,\ldots,u_d)$  is a {\em vector of orthogonality} for the $d$-OPS $\{P_n\}$. 
\end{defi}

For each $\{P_n\}$ being a $d$-OPS, the uniqueness of an associated vector of orthogonality is not guaranteed. In fact, given one of such vectors 
$(u_1,\ldots,u_d)$, $\{P_n\}$ is also a $d$-OPS with respect to 
$(v_1,\ldots,v_d)$ defined as
\begin{align*}
v_1&=u_1,\\
v_2&=u_2+\lambda_1^{(2)}u_1,\\
&\vdots\\
v_d&=u_d+\lambda^{(d)}_{d-1}u_{d-1}+\cdots+\lambda^{(d)}_{2}u_2+\lambda_1^{(d)}u_1,
\end{align*}
where $\lambda^{(j)}_{i}\in\mathbb{C},\,j=2,\ldots, d,\,i=1,\ldots , j-1,$ are randomly chosen. 

Conversely, if $(u_1,\ldots,u_d)$ is a vector of linear functionals and the sequence $\{P_n\}$ is a $d$-OPS with respect to $(u_1,\ldots,u_d)$, then $\{P_n\}$  is uniquely determined except a constant for each polynomial (see \cite{Maroni-1989}). In particular,  if $k_n$ denotes the leading coefficient of $P_n$ then the sequence $\{{\frac{1}{k_n}P_n}\}$ is called  sequence of monic $d$-orthogonal polynomials ({\em monic }$d$-OPS). Hereafter we alway deal with this kind of sequences of polynomials. 

A relevant characterization of $d$-OPS is the following (see also \cite{Maroni-1989}).

\begin{teo}[J. Van Iseghem~\cite{iseghem-1987}]
	 $\{P_n\},  n\in \mathbb{N},$  is a d-OPS if and only if there exist coefficients $\{a_{n,n-k},\,n\geq 0\,,\,k=0,\ldots,d\}$ such that the following $(d + 2)$-term recurrence relation is satisfied,
\begin{equation}
\label{reucrrence}
\left.
\begin{array}{rllcccc}
xP_{n}(x) & = & P_{n+1}(x)+\displaystyle\sum_{k=0}^{d}a_{n,n-k}P_{n-k}(x),\quad a_{n,n-d}\ne 0, \quad n\geq 0,\\
\\ P_{0}(x) & = & 1,\quad P_{-1}(x)=P_{-2}(x)=\cdots=P_{-d}(x)=0.	
\end{array}	
\right\}	 
\end{equation}
\end{teo}

In our research, starting from a $d$-OPS $\{P_n\}$ and an associated vector of orthogonality $(u_1,u_2,\ldots, u_d),$ we introduce $d$ new vectors of functionals $(u^{(m)}_1,u^{(m)}_2,\ldots,u^{(m)}_d),\,m=1,\ldots,d$, whose study is the object of this paper. In fact, for $m=1,$  it is possible to take a vector of functionals $(u^{(1)}_1,u^{(1)}_2,\ldots,u^{(1)}_d)$ verifying 
\begin{equation}\label{func(1)}
(x-a)u^{(1)}_1=u_d,\qquad u^{(1)}_i=u_{i-1},\quad i=2,\ldots, d.
\end{equation}
In this case
\begin{equation}
u_1^{(1)}=\dfrac{u_d}{x-a}+M_1\delta_a\,,
\label{*}
\end{equation}
where $\delta_a$ is the Dirac Delta function supported in $x=a$ (this is,
$
\langle\delta_a,q\rangle=q(a)
$
for any 
$q\in \mathcal{P}[x]$)
and $M_1\in \mathbb{C}$ is an arbitrary value. Therefore, $u^{(1)}_1$ is not uniquely determined from  $(u_1,u_2,\ldots, u_d)$. Similary, if $m=1,2,\ldots, d-1,$ and $(u^{(m)}_1,u^{(m)}_2,\ldots,u^{(m)}_d)$ is a vector of orthogonality for the $d$-OPS  $\{P_n^{(m)}\}$, then $(u^{(m+1)}_1,u^{(m+1)}_2,\ldots,u^{(m+1)}_d)$  verifying 
\begin{equation}\label{func(11)}
(x-a)u^{(m+1)}_1=u^{(m)}_d,\qquad u^{(m+1)}_i=u^{(m)}_{i-1},\quad i=2,\ldots, d.
\end{equation}
is a new vector of functionals. As in (\ref{*}) more than one vector of functionals can be defined, being
$$
u^{(m+1)}_1=\dfrac{u^{(m)}_d}{x-a}+M_{m+1}\delta_a
$$
with $M_2,\ldots, M_d\in \mathbb{C}$. 

From the above, the vectors of functionals $(u^{(r)}_1,u^{(r)}_2,\ldots,u^{(r)}_d)$ and $(u^{(r+q)}_1,u^{(r+q)}_2,\ldots,u^{(r+q)}_d)$, $r\in\{0,1,\ldots,d-1\},\,q\in\{1,\ldots,d\}$, are related by
\begin{equation} 
\label{nuevo6}
u_i^{(r)}=
\begin{cases}
u_{i+q}^{(r+q)},&i=1,\ldots, d-q.	\\[6pt]
(x-a)u_{i-d+q}^{(r+q)},&i=d-q+1,\ldots, d,
\end{cases}
\end{equation} 
The vectors of functionals  $(u^{(m)}_1,u^{(m)}_2,\ldots,u^{(m)}_d),\,m=1,\,\ldots, d,$ constructed in  (\ref{func(1)})-(\ref{func(11)}) extend the concept of Geronimus transformation of $(u_1,\ldots, u_d)$  introduced in \cite {Geronimus} and studied in several later works in the case $d=1$ (see  \cite{Bue1,DereM,Zhe-1997} for instance).

The values chosen for $M_1,\ldots, M_d$ can help to determine the regularity of the corresponding vectors of functionals. As we explain later, the vector of orthogonality $(u_1,u_2,\ldots, u_d)$  from which we start plays an essential role in the construction of $(u^{(m)}_1,u^{(m)}_2,\ldots,u^{(m)}_d),\,m=1,\,\ldots, d.$ One of our goals is to characterize the regularity of the new vectors of functionals verifying (\ref{func(1)})-(\ref{func(11)}), which is done in the following result.

\begin{teo}\label{representationP(m)-P}
	Let $(u_1,\ldots, u_d)$ be a  vector of orthogonality and let $\{P_n\}, n\in \mathbb{N},$ be its corresponding
	$d$-OPS. Then $(u^{(m)}_1,\ldots, u^{(m)}_{d}), $ $m=1,\ldots, d,$  is regular  if and only if  
	$d^{(m)}_n\neq 0,$  $n\in \mathbb{N}$, where
	
$$d^{(m)}_n=
\left\{
\begin{array}{ccc}
\left|\begin{array}{ccc}
	\langle{u_1^{(m)},P_{0}}\rangle&\cdots&\langle{u_n^{(
			m)},P_{0}}\rangle\\
	\vdots&\cdots&\vdots\\
	\langle{u_1^{(m)},P_{n-1}}\rangle&\cdots&\langle{u_n^{(m)},P_{n-1}}\rangle
	\end{array}\right|& ,&  m>n,\\\\
\left|\begin{array}{ccc}
	\langle{u_1^{(m)},P_{n-m}}\rangle&\cdots&\langle{u_m^{(
			m)},P_{n-m}}\rangle\\
		\vdots&\cdots&\vdots\\
	\langle{u_1^{(m)},P_{n-1}}\rangle&\cdots&\langle{u_m^{(m)},P_{n-1}}\rangle
	\end{array}
\right|	&, & m\leq n.
\end{array}
\right.
$$
	
In such  case the monic $d$-OPS, $\{P_n^{(m)}\}$, associated with  $(u^{(m)}_1,\ldots, u^{(m)}_{d}) $ is given by
	
\begin{equation}\label{cramer} 
	P_n^{(m)}(x)=\frac{1}{d^{(m)}_n}\begin{vmatrix}
	\langle{u_1^{(m)},P_{n-m}}\rangle&\cdots&\langle{u_m^{(m)},P_{n-m}}\rangle&P_{n-m}(x)\\[5pt]
	\vdots&&\vdots&\vdots\\
	\langle{u_1^{(m)},P_{n-1}}\rangle&\cdots&\langle{u_m^{(m)},P_{n-1}}\rangle&P_{n-1}(x)\\[5pt]
	\langle{u_1^{(m)},P_{n}}\rangle&\cdots&\langle{u_m^{(m)},P_{n}}\rangle&P_{n}(x)\\
	\end{vmatrix}, \quad n\geq m.
	\end{equation}
\end{teo}

Theorem \ref{representationP(m)-P} extends previous results corresponding to the case $d=1$. More specifically, the following corollary is an immediate consequence, which was proved in \cite{Dere} with other arguments. 

\begin{coro}
Let $(u_1,\ldots, u_d)$ be a regular vector of  linear functionals and consider $\{P_n\},\ n\in \mathbb{N},$ its corresponding
$d$-OPS. Then $(u^{(1)}_1,\ldots ,u^{(1)}_{d})$  is  regular  if only if    $\langle u^{(1)}_1,P_{n} \rangle\ne 0$  for each $n\geq 0.$ In such  case, $\{P_n^{(1)}\}$ exists and 

$$P_n^{(1)}(x)=\frac{1}{\langle u^{(1)}_1,P_{n-1}\rangle}\begin{vmatrix}
P_n(x)&\langle u^{(1)}_1,P_{n}\rangle\\[10pt]
P_{n-1}(x)&\langle u^{(1)}_1,P_{n-1}\rangle
\end{vmatrix}, \quad n\in \mathbb{N}.$$
\end{coro}

The $(d+2)$-recurrence relation \eqref{reucrrence} motives our interes in the following Hessenberg banded matrix, univocally determined by each $d$-OPS. Set
$$
J=
\left(
\begin{array}{cccccc}
a_{0,0}&1        &      &             &\\
a_{1,0}&a_{1,1}  &1     &             & \\
\vdots & \vdots  &\ddots&\ddots       &\\
a_{d,0}&a_{d,1}  &\cdots&a_{d,d}      &1\\
0      &a_{d+1,1}&      &&\ddots&\ddots\\
   &\ddots&\ddots\\
\end{array}
\right),	 \qquad a_{d+i,i}\ne 0, \quad i=0,1,\ldots,
$$
whose entries are the coefficients of the recurrence relation \eqref{reucrrence}. Assume that $(u_1,\ldots, u_d)$ is a vector of orthogonality  associated with $J$ and consider that $(u^{(m)}_1,\ldots, u^{(m)}_d),\,m=1,\ldots, d,$ are vectors of orthogonality verifying  (\ref{func(1)})-(\ref{func(11)}). Then, for each $m=1,\ldots, d,$ there exists a $d$-OPS $\{P^{(m)}_n\}$ verifying the corresponding $(d+2)$-term recurrence relation whose coefficients define a new Hessenberg banded matrix
$$
J^{(m)}=
\left(
\begin{array}{cccccc}
a^{(m)}_{0,0}&1        &      &             &\\
a^{(m)}_{1,0}&a^{(m)}_{1,1}  &1     &             & \\
\vdots & \vdots  &\ddots&\ddots       &\\
a^{(m)}_{d,0}&a^{(m)}_{d,1}  &\cdots&a^{(m)}_{d,d}      &1\\
0      &a^{(m)}_{d+1,1}&      &&\ddots&\ddots\\
   &\ddots&\ddots\\
\end{array}
\right),	 \qquad a^{(m)}_{d+i,i}\ne 0, \quad i=0,1,\ldots,
$$
associated with $(u^{(m)}_1,\ldots, u^{(m)}_d)$. We are concerned to finding relations between $J$ and $J^{(m)}, \,m=1,\ldots, d$. Our second main result is the following.
\begin{teo}
\label{principal}
 Let $\left(u^{(m)}_1,u^{(m)}_2,\ldots,u^{(m)}_d\right),\,m=0,1,\ldots, d,$ be a set of vector of orthogonality verifying  (\ref{func(1)})-(\ref{func(11)}) and assume that $\{P^{(m)}_n\},n\in\mathbb{N}\,,m=0,1,\ldots, d,$ are the corresponding $d$-OPS (where  
 $(u^{(0)}_1,u^{(0)}_2,\ldots,u^{(0)}_d)= (u_1,u_2,\ldots,u_d)$ and $\{P_n^{(0)}\}=\{P_n\}$). Then there exist $d$ lower triangular matrices
\begin{equation}
 \label{Ls}
L^{(m)}=
\left(
\begin{array}{cccc}
1 \\
\gamma^{(m)}_{1} & 1 \\
 &  \gamma^{(m)}_{2} & \ddots  \\
 & & \ddots
\end{array}
\right),\quad m=1,\ldots, d,
\end{equation}
and there exists an upper triangular matrix 
$$
U=
\left(
\begin{array}{cccc}
u_1& 1 \\
& u_2 & 1 \\
 & & \ddots  & \ddots  \\
\end{array}
\right)
$$
 such that
\begin{equation}
\label{otrooctavo}
J^{(m)}-aI=L^{(m)}\cdots L^{(1)}UL^{(d)}\cdots L^{(m+1)},\qquad m=1,\ldots, d.
\end{equation}
 
\end{teo}

In the rest of the paper we consider a regular vector of functionals $(u_1,\ldots, u_{d}) $ and we assume that $\{P_n\}$ is the monic $d$-OPS associated with $(u_1,\ldots, u_{d}) $. 

In Section \ref{section2} some connections between the $d$ sequences of polynomials $\{P^{(m)}_n\}$ are studied, as well as some factorizations for matrices $J^{(m)}-aI$, in the case that the vectors of functionals $(u^{(m)}_1,\ldots, u^{(m)}_{d}), m=1,\ldots, d, $ are regular. The auxiliary results introduced in Section \ref{section2} are the basis of the proof of Theorem \ref{principal}, which will be easily obtained from these relationships. Finally, Section \ref{section3} is devoted to prove Theorem \ref{representationP(m)-P} and  Theorem \ref{principal}.

\section{Geronimus transformation for vectors of orthogonality}\label{section2}

Theorem \ref{representationP(m)-P} provides some conditions easy to check which guarantee the regularity of the analized vectors of functionals. For instance, if $m=1$ we have that $(u^{(1)}_1,\ldots, u^{(1)}_{d}), $ is regular when 
$$
d^{(1)}_n=\langle u^{(1)}_1, P_{n-1}\rangle \neq 0\,,\quad n\in \mathbb{N}.
$$
This is, from \eqref{*}, 
\begin{equation*}
\label{nuevo*}
\langle \frac{u_d}{x-a},P_{n-1}\rangle +M_1 P_{n-1}(a)\neq 0\,,\quad n\geq 1.
\end{equation*}
Hence, condition
\begin{equation*}
\label{condition}
P_n(a)\neq 0\,,\quad n\in \mathbb{N},
\end{equation*}
permits to choose $M_1\in \mathbb{C}$,
\begin{equation}
\label{elegirM1}
M_1\neq -\frac{\langle\displaystyle \frac{u_d}{x-a},P_{n-1}\rangle }{P_n(a)}\,,\quad \forall n\in \mathbb{N},
\end{equation}
for obtaining a vector of orthogonality $(u^{(1)}_1,\ldots, u^{(1)}_{d})$. 

We can iterate the procedure, obtaining the vectors of orthogonality $(u^{(m)}_1,\ldots, u^{(m)}_{d})$ and their corresponding $d$-OPS, $\{P^{(m)}_n\}$, for $m=1,2,\ldots, k,\,k\leq d,$ when 
\begin{equation}
\label{a-condicion}
P_n(a),\,P^{(1)}_n(a),\,\ldots,\,P^{(m-1)}_n(a)\neq 0\,,\quad \forall n\in \mathbb{N}.
\end{equation}
Nevertheless, \eqref{a-condicion} is not a necessary condition for constructing $P^{(k)}_n$. Indeed, $(u^{(k)}_1,u^{(k)}_2,\ldots,u^{(k)}_d)$ can be directly defined from $(u_1,u_2,\ldots,u_d)$ verifying \eqref{nuevo6} even if the vectors of functionals $(u^{(m)}_1,u^{(m)}_2,\ldots,u^{(m)}_d),\,m=1,\ldots, k-1,$ are not previously defined. 

On the other hand, we emphasize that the vectors of functionals verifying (\ref{func(1)})-(\ref{func(11)}) are related by \eqref{nuevo6}, which is independent on whether  $(u^{(r)}_1,\ldots, u^{(r)}_{d})$ is or not regular. 

In the remainder of this section we take $r\in \{0,\ldots, d-1\}$  and $q\in\{1,\ldots, d-r\}$ fixed and we assume that $\left(u^{(k)}_1,u^{(k)}_2,\ldots,u^{(k)}_d\right),\,k=r,r+q, $ are regular. Then for $k=r,r+q,$ there exists $\{P^{(k)}_n\}, n\in \mathbb{N},$ which is a $d$-OPS with respect to $(u^{(k)}_1,u^{(k)}_2,\ldots,u^{(k)}_d)$. The sequence $\{P^{(k)}_n\}$ satisfy a $(d+2)$-term recurrence relation,
\begin{equation}
\label{recurrence}
\left.
\begin{array}{cclcc}
xP^{(k)}_{n}(x)&= & P^{(k)}_{n+1}(x)+\displaystyle\sum_{s=0}^{d}a^{(k)}_{n,n-s}P^{(k)}_{n-s}(x),\quad a^{(k)}_{n,n-d}\ne 0, \quad n\geq 0,\\\\
P^{(k)}_{0}(x)& = & 1, \quad P^{(k)}_{-1}(x)=P^{(k)}_{-2}(x)=\cdots=P^{(k)}_{-d}(x)=0.
\end{array}	
\right\}
\end{equation}
Next, we analyze the relation between the sequences $\{P^{(r)}_n\}$  and $\{P^{(r+q)}_n\}$. (We recall that we understand $P^{(0)}_n=P_n$ and $u_i^{(0)}=u_i,\,i=1,\ldots, d.$)

\begin{lem}
\label{lema1} With the above notation, for each $n\in \mathbb{N}$ there exists  a set of complex numbers  $\{\gamma_{n,k}^{(r,q)}: \,k=n-q,\ldots, n-1,\, \gamma_{n,n-q}^{(r,q)}\}~\neq~0$ such that 
 \begin{equation}\label{rec}
	 P^{(r+q)}_{n}(x)=P^{(r)}_{n}(x)+\sum_{s=n-q}^{n-1}{\gamma^{(r,q)}_{n,s}}P^{(r)}_{s}(x),\quad  n\geq 0.
	 \end{equation}
\end{lem}

\noindent
\underline{Proof.- }
From \eqref{nuevo6}, for $j\in \{1,\ldots,d\}$ fixed we have 
\begin{equation}
\label{10}
\left.
\begin{array}{ccl}
\langle u^{(r)}_i,\,x^kP^{(r+q)}_{dn+j}\rangle=0&,& k=0,1,\ldots,n-1 \\\\
\langle u^{(r)}_i,\,x^{n}P^{(r+q)}_{dn+j}\rangle=0&,& i+q\leq j
\end{array}
\right\}\,,\quad i=1,\ldots, d-q,
\end{equation}
and 
\begin{equation}
\label{11}
\left.
\begin{array}{lcl}
\langle u^{(r)}_i,\,x^kP^{(r+q)}_{dn+j}\rangle=0&,& k=0,1\ldots, n-2 \\\\
\langle u^{(r)}_i,\,x^{n-1}P^{(r+q)}_{dn+j}\rangle=0&,& i-d+q\leq j
\end{array}
\right\}\,,\quad i=d-q+1,\ldots, d.
\end{equation}

On the other hand, for each $n\in \mathbb{N}$ there exists a set $\{\gamma^{(r,q)}_{dn+j,s}\,:\,s=0,\ldots, dn+j-1 \}$ such that 
\begin{equation}
\label{12}
 P^{(r+q)}_{dn+j}(x)=P^{(r)}_{dn+j}(x)+\sum_{s=0}^{dn+j-1}{\gamma^{(r,q)}_{dn+j,s}}P^{(r)}_{s}(x)\,.
 \end{equation}
Thus, using  \eqref{10} and \eqref{11}, we have 

\begin{equation} 
\label{nuevo15}
\left.\begin{array}{ccc}
\displaystyle\sum_{s=0}^{kd}\gamma^{(r,q)}_{dn+j,s}\ \langle u_1^{(r)},x^kP^{(r)}_s\rangle&=&0,\\
	\displaystyle\sum_{s=0}^{kd+1}\gamma^{(r,q)}_{dn+j,s}\langle u_2^{(r)},x^kP^{(r)}_s\rangle&=&0,\\
	\vdots && \vdots\\
	\displaystyle\sum_{s=0}^{(k+1)d-1}\gamma^{(r,q)}_{dn+j,s}\langle u_d^{(r)},x^kP^{(r)}_s\rangle&=&0
\end{array}\right\},\quad  k=0,\ldots,n-2.
\end{equation}

Due to $\langle u^{(r)}_i, \,x^{k }P^{(r)}_{s}\rangle\neq 0$ for $s=dk+i-1$, taking successively $k=0,1,\ldots, n-2$ in \eqref{nuevo15} we obtain 
\begin{equation}
\label{nuevo18}
\gamma^{(r,q)}_{dn+j,0}=\gamma^{(r,q)}_{dn+j,1}=\cdots=\gamma^{(r,q)}_{dn+j,d(n-1)-1}=0.
\end{equation}
From this and the first condition of \eqref{10}
\begin{equation*}
\label{otraotra12}
\left.
\begin{array}{cll}
\displaystyle\sum_{s=d(n-1)}^{d(n-1)}{\gamma^{(r,q)}_{dn+j,s}}\langle u^{(r)}_1, \,x^{n-1} P^{(r)}_{s}\rangle & = &0 \\
\displaystyle\sum_{s=d(n-1)}^{d(n-1)+1}{\gamma^{(r,q)}_{dn+j,s}}\langle u^{(r)}_2, \,x^{n-1} P^{(r)}_{s}\rangle & = &0 \\
\vdots && \vdots  \\
\displaystyle\sum_{s=d(n-1)}^{d(n-1)+d-q-1}{\gamma^{(r,q)}_{dn+j,s}}\langle u^{(r)}_{d-q}, \,x^{n-1} P^{(r)}_{s}\rangle & = &0 \\
\end{array}
\right\}.
\end{equation*}
Then,  
\begin{equation}
\label{otronuevo18}
\gamma^{(r,q)}_{dn+j,d(n-1)}=\gamma^{(r,q)}_{dn+j,d(n-1)+1}=\cdots=\gamma^{(r,q)}_{dn+j,dn-q-1}=0\,.
\end{equation}
Using \eqref{nuevo18} and \eqref{otronuevo18}, \eqref{12} reduces to 
\begin{equation}
 P^{(r+q)}_{dn+j}(x)=P^{(r)}_{dn+j}(x)+\sum_{s=dn-q}^{dn+j-1}{\gamma^{(r,q)}_{dn+j,s}}P^{(r)}_{s}(x)\,,\quad j=1,\ldots,d\,.
\label{otro19}
\end{equation}
To obtain \eqref{rec}, firstly assume $j\leq q$. Then from \eqref{11},
 $$
\left.
\begin{array}{cccclll}
\gamma^{(r,q)}_{dn+j,dn-q}\langle u^{(r)}_{d-q+1}, \,x^{n-1} P^{(r)}_{dn-q}\rangle &&& = &0 \\\\
\gamma^{(r,q)}_{dn+j,dn-q}\langle u^{(r)}_{d-q+2}, \,x^{n-1} P^{(r)}_{dn-q}\rangle &+&\gamma^{(m)}_{dn+j,dn-q+1}\langle u^{(r)}_d-q+2, \,x^{n-1} P^{(r)}_{dn-q+1}\rangle & = &0 \\
\vdots &\vdots &\vdots  \\
\gamma^{(r,q)}_{dn+j,dn-q}\langle u^{(r)}_{d-q+j}, \,x^{n-1} P^{(r)}_{dn-q}\rangle&+&\cdots+\gamma^{(r,q)}_{dn+j,dn+j-q-1}\langle u^{(r)}_{d-q+j}, \,x^{n-1} P^{(r)}_{dn+j-q-1}\rangle &= &0 \\
\end{array}
\right\}
$$
Thus
\begin{equation}
\label{nuevo19}
\gamma^{(r,q)}_{dn+j,dn-q}=\gamma^{(r,q)}_{dn+j,dn-q+1}=\cdots =\gamma^{(r,q)}_{dn+j,dn+j-q-1}=0.
\end{equation}
Moreover, since   $(x-a)u^{(r+q)}_{j+1}=u^{(r)}_{d-q+j+1},$ 
$$\gamma^{(r,q)}_{dn+j,dn+j-q}=\dfrac{\langle u^{(r)}_{d-q+j+1},x^{n-1}P^{(r+q)}_{dn+j}\rangle}{\langle u^{(r)}_{d-q+j+1},x^{n-1}P^{(r)}_{dn+j-q}\rangle}=\dfrac{\langle u^{(r+q)}_{j+1},(x-a)x^{n-1}P^{(r+q)}_{dn+j}\rangle}{\langle u^{(r)}_{d-q+j+1},x^{n-1}P^{(r)}_{dn+j-q}\rangle}\ne 0.$$

Secondly, assume $j>q$. Then, from \eqref{10} and \eqref{11},
$$
\left.
\begin{array}{cccclllllll}
\gamma^{(r,q)}_{dn+j,dn-q}\langle u^{(r)}_{d-q+1}, \,x^{n-1} P^{(r)}_{dn-q}\rangle&&&&&=0 \\
\vdots  \\
\gamma^{(r,q)}_{dn+j,dn-q}\langle u^{(r)}_{d}, \,x^{n-1} P^{(r)}_{dn-q}\rangle&+&\cdots &+&\gamma^{(r,q)}_{dn+j,dn-1}\langle u^{(r)}_{d}, \,x^{n-1} P^{(r)}_{dn-1}\rangle & = 0 \\\\
\gamma^{(r,q)}_{dn+j,dn-q}\langle u^{(r)}_{1}, \,x^{n} P^{(r)}_{dn-q}\rangle&+&\cdots&+&\gamma^{(r,q)}_{dn+j,dn}\langle u^{(r)}_{1}, \,x^{n} P^{(r)}_{dn}\rangle & =0 \\
\vdots &\vdots& & \vdots &\qquad \qquad\vdots \\
\gamma^{(r,q)}_{dn+j,dn-q}\langle u^{(r)}_{j-q}, \,x^{n} P^{(r)}_{dn-q}\rangle&+&\cdots&+&\gamma^{(r,q)}_{dn+j,dn+j-q-1}\langle u^{(r)}_{j-q}, \,x^{n} P^{(r)}_{dn+j-q-1}\rangle &= 0 
\end{array}
\right\}
$$
Again, this implies
\begin{equation}
\label{otronuevo19}
\gamma^{(r,q)}_{dn+j,dn-q}=\gamma^{(r,q)}_{dn+j,dn-q+1}=\cdots =\gamma^{(r,q)}_{dn+j,dn+j-q-1}=0.
\end{equation}
Thus, taking \eqref{nuevo19}-\eqref{otronuevo19} into consideration, from \eqref{otro19} we arrive to \eqref{rec}.

Finally, taking into account that $u^{(r+q)}_{j+1}=u^{(r)}_{j-q+1}$ (see \eqref{nuevo6}),

 $$\gamma^{(r,q)}_{dn+j,dn+j-q-1}=\dfrac{\langle u^{(r)}_{j-q+1},x^nP^{(r+q)}_{dn+j}\rangle}{\langle u^{(r)}_{j-q+1},x^nP^{(r)}_{dn+j-q}\rangle}=\dfrac{\langle u^{(r+q)}_{j+1},x^nP^{(r+q)}_{dn+j}\rangle}{\langle u^{(r)}_{j-q+1},x^nP^{(r)}_{dn+j-q}\rangle}\ne 0$$
and the result is proved.   $\hfill \square$

Lemma \ref{lema1} provides the expression of the polynomials $\{ P^{(r+q)}_n\}$ in terms of $\{ P^{(r)}_n\}$. Reciprocally, the following auxiliary result provides each polynomial $P^{(r)}_n$ in terms of the sequence  $\{ P^{(r+q)}_n\}$.

\begin{lem}
\label{lema2}
For each $n\in \mathbb{N}$ we have
\begin{equation}
\label{13}
(x-a) P^{(r)}_{n}(x)=P^{(r+q)}_{n+1}(x)+\sum_{s=n-d+q}^n \alpha^{(r,q)}_{n,s}P^{(r+q)}_{s}(x)\,,
\end{equation}
where $\alpha^{(r,q)}_{n,n-d+q}\neq 0$.
\end{lem}

\noindent
\underline{Proof.- }
For $n\in \mathbb{N}$ and $j\in\{1,\ldots, d\}$ the polynomial $(x-a) P^{(r)}_{dn+j}(x)$ can be obtained in terms of $\{P^{(r+q)}_{n}\}$ as 
$$
(x-a) P^{(r)}_{dn+j}(x)=P^{(r+q)}_{dn+j+1}(x)+\sum_{k=0}^{dn+j} \alpha^{(r,q)}_{n,k}P^{(r+q)}_{k}(x)\,.
$$
From \eqref{nuevo6} we have 
\begin{equation}
\label{101}
\left.
\begin{array}{ccl}
\langle u^{(r+q)}_i,\,(x-a)x^kP^{(r)}_{dn+j}\rangle=0&,& k=0,1\ldots n-1 \\\\
\langle u^{(r+q)}_i,\,(x-a)x^{n}P^{(r)}_{dn+j}\rangle=0&,& d-q+i\leq j
\end{array}
\right\}\,,\quad  i=1,\ldots, q,
\end{equation}
and 
\begin{equation}
\label{111}
\langle u^{(r+q)}_i,\,x^kP^{(r)}_{dn+j}\rangle=0,\quad  k=0,1\ldots, n-1 
,\quad i=q+1,\ldots, d,
\end{equation}
Then, similar arguments to those used in the proof of Lemma \ref{lema1} lead to (\ref{13}). 
$\hfill \square$

\begin{rem}
 \eqref{rec} can be rewritten as 
 \begin{equation}
 \label{nota1}
 {\bf P}^{(r+q)}(x)=L^{(r,q)}{\bf P}^{(r)}(x),
 \end{equation}
 where 
 $${\bf P}^{(k)}(x)=\left(  P^{(k)}_{0}(x),  P^{(k)}_{1}(x), \ldots\right)^T,\,k=r,r+q,$$ 
 and $ L^{(r,q)}$ is a lower triangular banded matrix,
 \begin{equation*}
	L^{(r,q)}=\left(
	\begin{matrix}
	1 \\
		\gamma^{(r,q)}_{1,0}&1        &      &             &\\[5pt]
		\vdots & \vdots  &\ddots&\\
		\gamma^{(r,q)}_{q,0}&\gamma^{(r,q)}_{q,1}  &\cdots&1\\
		0      &\gamma^{(r,q)}_{q+1,1}&   \ddots   &&\ddots&\\
		 &0  &&\ddots&&\\
		  &  &\ddots&\ddots&\ddots&
	\end{matrix}
	\right), 
\end{equation*}

In the case that $(u_1^{(k)},\ldots ,u_{d}^{(k)}),\, k\in\{r+1,\ldots,r+q-1\},$ are regular vectors of functionals then the same above reasoning can be done replacing $q$ by 1 and substituting $r$ succesively by $r+1, \ldots, r+q-1$  in \eqref{nota1}. In this way,  
$$
\begin{array}{ccccllll}
 {\bf P}^{(r+2)}(x)& = & L^{(r+1,1)}{\bf P}^{(r+1)}(x)&=&L^{(r+1,1)}L^{(r,1)}{\bf P}^{(r)}(x)\\
 {\bf P}^{(r+3)}(x)& = & L^{(r+2,1)}{\bf P}^{(r+2)}(x)&=&L^{(r+2,1)}L^{(r+1,1)}L^{(r,1)}{\bf P}^{(r)}(x)\\
\vdots && \vdots & &\qquad\vdots \\
{\bf P}^{(r+q)}(x)& = & L^{(r+q-1,1)}{\bf P}^{(r+q-1)}(x)&=& L^{(r+q-1,1)}\cdots L^{(r,1)}{\bf P}^{(r)}(x).
\end{array}
$$
Hence,
\begin{equation}
\label{nuevo22}
{\bf P}^{(r+q)}(x)=L^{(r+q)} L^{(r+q-1)}\cdots L^{(r+1)}{\bf P}^{(r)}(x).
\end{equation}
In \eqref{nuevo22} and in the sequel we denote 
\begin{equation}
\label{otronuevo22}
L^{(s+1)}:= L^{(s,1)}\,,\quad s=r,r+1,\ldots, r+q-1.
\end{equation}
Notice that these matrices $L^{(s+1)}$ are bi-diagonal as in (\ref{Ls}), where we are writing 
$
\gamma_i^{(s+1)}:=\gamma_{i+1,i}^{(s,1)},\,i=1,2\ldots
$
 (We remember that we are assuming that 
 $(u_1^{(r)},\ldots ,u_{d}^{(r)})$  
 and
  $(u_1^{(r+q)},\ldots, u_{d}^{(r+q)})$ 
are regular vectors of functionals). 
We underline that, under these conditions, the following factorization of $L^{(r,q)}$ is verified,
\begin{equation}
\label{nuevo23}
L^{(r,q)}=L^{(r+q)} L^{(r+q-1)}\cdots L^{(r+1)}.
\end{equation}

In the same way, (\ref{13}) can be rewritten as 
\begin{equation}
 \label{nota2}
(x-a){\bf P}^{(r)}(x)=N^{(r,q)}{\bf P}^{(r+q)}(x)\,,
\end{equation}
where $N^{(r,q)}$ is a lower Hessenberg $(d+2-q)$-banded matrix,

\begin{equation}
\label{nuevo21}
N^{(r,q)}=
\left(
	\begin{matrix}
		\alpha^{(r,q)}_{0,0}&1        &      &             &\\[5pt]
		\alpha^{(r,q)}_{1,0}&\alpha^{(r,q)}_{1,1}  &1     &             & \\
		\vdots & \vdots  &\ddots&\ddots       &\\
		\alpha^{(r,q)}_{d-q,0}&\alpha^{(r,q)}_{d-q,1}  &\cdots&\alpha^{(r,q)}_{d-q,d-q}      &1\\
		0      &\alpha^{(r,q)}_{d-q+1,1}&     &\ddots&&\ddots\\
		 &0  & \ddots&&\ddots\\
		  &  &\ddots&\ddots&&\ddots
	\end{matrix}
	\right), 
	\quad \alpha^{(m)}_{d-q+i,i}\ne 0\,,	i=0,1,\ldots
\end{equation}

\end{rem}

Lemma \ref{lema1} and Lemma \ref{lema2} provide relationships between the matrices $J^{(r)}$ and $J^{(r+q)}$. This is sumarized in the following result. 

\begin{teo}\label{teorema4}
	With the above notation we have
	\begin{equation}\label{J}
	J^{(r)}-aI=N^{(r,q)}L^{(r,q)}
	\end{equation}
and
	\begin{equation}\label{jm}
	J^{(r+q)}-aI=L^{(r,q)}N^{(r,q)}.
	\end{equation}
\end{teo}

\noindent
\underline{Proof.- }
As we have done
	\begin{equation*}
	(J^{(r)}-aI){\bf P}^{(r)}(x)=(x-a){\bf P}^{(r)}(x)=N^{(r,q)} {\bf P}^{(r+q)}(x)=N^{(r,q)}L^{(r,q)}{\bf P}^{(r)}(x).
	\end{equation*}
	From here, and taking into account that  $\{P^{(r)}_n\}$ is a basis of $\mathcal{P}[x]$, \eqref{J} is verified. 
	
	To prove \eqref{jm} notice that
	
	\begin{equation*}
	( J^{(r+q)}-aI){\bf P}^{(r+q)}(x)=(x-a)L^{(r,q)}{\bf P}^{(r)}(x)=L^{(r,q)}N^{(r,q)}{\bf P}^{(r+q)}(x).
	\end{equation*}
Using again the fact that $\{P_n^{(r+q)}\}, n\in \mathbb{N},$ is a basis of $\mathcal{P}[x]$, \eqref{jm} follows. 
$\hfill\square$

\section{Proofs of the main results}\label{section3}
	
\subsection{Proof of Theorem \ref{representationP(m)-P}}

 We will give the proof when $n\geq m.$ The case $n<m$ follows in a similar way. 
 
 In the first place, we take $m\in \{1,\ldots, d\}$ and we suppose that $(u_1^{(m)},\ldots,u_{d}^{(m)})$ is regular. From Lemma \ref{lema1} for $r=0$ and $q=m$, we have
$$
{\bf P}^{(m)}(x)=L^{(0,m)}{\bf P}(x)
$$
(see \eqref{nota1}), where $L^{(0,m)}$ is a lower triangular $(m+1)$-banded matrix as in \eqref{nota2}. 

This is, 
 \begin{equation}\label{123}
 P_n^{(m)}(x)=P_{n}(x)+\gamma^{(0,m)}_{n,n-1}P_{n-1}(x)+\cdots+\gamma^{(0,m)}_{n,n-m}P_{n-m}(x)
 \end{equation}
(see \eqref{rec}).
 
 We recall that  $
\gamma^{(0,m)}_{n,n-m}\neq 0
 $. Taking into account 
 $$\langle{u^{(m)}_k,P_n^{(m)}}\rangle=0,\quad  k=1,\ldots, m,$$ 
 and \eqref{123}, we have 
\begin{equation}
\label{33.1}
\left(
 \begin{matrix}
 \langle u_1^{(m)},P_{n-m}\rangle&\cdots&\langle u_1^{(m)},P_{n-1}\rangle\\
 \vdots&\cdots&\vdots\\
 \langle{u_m^{(m)},P_{n-m}}\rangle&\cdots&\langle u_m^{(m)},P_{n-1}\rangle
 \end{matrix}
 \right)
 \left(
\begin{matrix}
 \gamma^{(0,m)}_{n,n-m}\\
 \vdots\\
  \gamma^{(0,m)}_{n,n-1}
 \end{matrix}
 \right)=
  -\left(
  \begin{matrix}
\langle{u_1^{(m)},P_n}\rangle\\
 \vdots\\
\langle u_m^{(m)},P_n\rangle
 \end{matrix}
 \right)\,.
\end{equation}
This is, 
$
(t_1,\ldots, t_m)=( \gamma^{(0,m)}_{n,n-m},\ldots,  \gamma^{(0,m)}_{n,n-1})
$
is a solution of the linear system
\begin{equation}
\label{sistema}
\left.
 \begin{array}{cccccccc}
 \langle u_1^{(m)},P_{n-m}\rangle t_1 & + & \cdots& + &\langle u_1^{(m)},P_{n-1}\rangle t_m & = &  -\langle u_1^{(m)},P_{n}\rangle \\
 \vdots& \vdots&& \vdots&\vdots& \vdots&\vdots\\
  \langle u_m^{(m)},P_{n-m}\rangle t_1 & + & \cdots& + &\langle u_m^{(m)},P_{n-1}\rangle t_m & = &  -\langle u_m^{(m)},P_{n}\rangle
 \end{array}
 \right\}
\end{equation}
We check that \eqref{sistema} has a unique solution. In fact, if $\left( t^{(m)}_{n},\ldots,t^{(m)}_{n-m+1}\right)^T$ is another solution then we define the polynomial
$$Q_n^{(m)}(x)=P_{n}(x)+t^{(m)}_{n}P_{n-1}(x)+\cdots+t^{(m)}_{n-m+1}P_{n-m}(x).$$ 
In what follows we will prove
\begin{equation}
\label{333.1}
\langle u_k^{(m)},(x-a)^rQ_n^{(m)}\rangle=~0,\quad k=1,\ldots, d,\quad n\geq dr+k, \quad n\geq m,\quad r\geq 0\,.
\end{equation}
In the case that $k\in\{1,\ldots,m\}$, then $\langle u_k^{(m)},Q_n^{(m)}\rangle=~0$ because 
	$\left( t^{(m)}_{n},\ldots,t^{(m)}_{n-m+1}\right)^T$
	is a solution of \eqref{33.1}. Therefore, for $r>0$ and $n\geq dr+k$, since \eqref{nuevo6} we obtain
\begin{eqnarray}	
\label{otro35}
\langle u_k^{(m)},(x-a)^rQ_n^{(m)}\rangle  & = & \langle(x-a)u_k^{(m)},(x-a)^{r-1}P_n\rangle+ t^{(m)}_{n}\langle(x-a)u_k^{(m)},(x-a)^{r-1}P_{n-1}\rangle\nonumber \\
&+ & \cdots+t^{(m)}_{n-m+1}\langle(x-a)u_k^{(m)},(x-a)^{r-1}P_{n-m}\rangle\nonumber\\
	& = & \langle u_{d-m+k},(x-a)^{r-1}P_n\rangle+t^{(m)}_{n}\langle u_{d-m+k},(x-a)^{r-1}P_{n-1}\rangle\nonumber\\
	& + &\cdots
	+t^{(m)}_{n-m+1}\langle u_{d-m+k},(x-a)^{r-1}P_{n-m}\rangle=0.
\end{eqnarray}		

On the contrary, if $k=m+p$  with $p\in \{1, \ldots , d-m\}$ and $n\geq dr+k~$ then, again since \eqref{nuevo6}, we obtain
\begin{eqnarray}
\label{otronuevo35}
\langle u_{k}^{(m)},(x-a)^rQ_n^{(m)}\rangle& =& 
\langle u_p,(x-a)^{r}P_n\rangle+t^{(m)}_{n}\langle u_p,(x-a)^{r}P_{n-1}\rangle\nonumber\\
& + & \cdots+t^{(m)}_{n-m+1}\langle u_p,(x-a)^{r}P_{n-m}\rangle=0. 
\end{eqnarray}
Then, since \eqref{otro35}-\eqref{otronuevo35} we see that \eqref{333.1} is verified for any $k\in \{1,\ldots, d\}$.

Obviously, \eqref{333.1}  implies $\langle u^{(m)}_k,\,x^rQ^{(m)}_n\rangle =0$, $n\geq rd+k.$
Doing a similar analysis, we obtain that for $n=dr+k-1\geq m$ is 
$\langle u^{(m)}_k,x^rQ^{(m)}_{dr+k-1}\rangle \ne 0.$
Thus, $\{Q_n^{(m)}\}, \,n\in \mathbb{N},$ is also a monic $d$-OPS, which contradicts the uniqueness of $\{P_n^{(m)}\}.$ 
This proves the uniqueness of solutions for the system \eqref{33.1}. Finally, applying the well-known Cramer Rule to solve \eqref{sistema} we arrive to \eqref{cramer}. 

Reciprocally, assume $d^{(m)}_n\ne 0$ and define $P^{(m)}_n$ as in \eqref{cramer}. We want to show that $\{P^{(m)}_n\}$ is the monic $d$-OPS with respect to $(u^{(m)}_1,\ldots,u^{(m)}_{d})$.

Take $r\in\mathbb{N}, \,k\in\{1,\ldots, d\}$. Then, for $n\geq dr+k$ we have

\begin{equation}
\label{(a)}
\langle{u_k^{(m)},(x-a)^rP_n^{(m)}}\rangle=\dfrac{1}{d^{(m)}_n}\begin{vmatrix}
\langle{u_1^{(m)},P_{n-m}}\rangle&\cdots&\langle{u_m^{(m)},P_{n-m}}\rangle&\langle{u_k^{(m)},(x-a)^rP_{n-m}}\rangle\\[5pt]
\vdots&&\vdots&\vdots\\
\langle{u_1^{(m)},P_{n-1}}\rangle&\cdots&\langle{u_m^{(m)},P_{n-1}}\rangle&\langle{u_k^{(m)},(x-a)^rP_{n-1}}\rangle\\[5pt]
\langle{u_1^{(m)},P_{n}}\rangle&\cdots&\langle{u_m^{(m)},P_{n}}\rangle&\langle{u_k^{(m)},(x-a)^rP_{n}}\rangle\\
\end{vmatrix}\,.
\end{equation}

If $k\leq m$ and $r=0,$ then we see on the right hande side of \eqref{(a)} that $\langle{u_k^{(m)},P_n^{(m)}}\rangle=0$ because two columns of the determinant coincide. If $k\leq m$ and $r\ne 0,$ then 

\begin{equation}
\label{(b)}
\langle{u_k^{(m)},(x-a)^rP_n^{(m)}}\rangle=\dfrac{1}{d^{(m)}_n}\begin{vmatrix}
\langle{u_1^{(m)},P_{n-m}}\rangle&\cdots&\langle{u_m^{(m)},P_{n-m}}\rangle&\langle{u_{d-m+k},(x-a)^{r-1}P_{n-m}}\rangle\\[5pt]
\vdots&&\vdots&\vdots\\
\langle{u_1^{(m)},P_{n-1}}\rangle&\cdots&\langle{u_m^{(m)},P_{n-1}}\rangle&\langle{u_{d-m+k},(x-a)^{r-1}P_{n-1}}\rangle\\[5pt]
\langle{u_1^{(m)},P_{n}}\rangle&\cdots&\langle{u_m^{(m)},P_{n}}\rangle&\langle{u_{d-m+k},(x-a)^{r-1}P_{n}}\rangle\\
\end{vmatrix}=0
\end{equation}
because the entries in the last column of the determinant are zero. 

If $k>m,$ then $k=m+p$ with $p\in \{1,\ldots,d-m\}$. Therefore, 
\begin{equation}
\label{(c)}
\langle{u_k^{(m)},(x-a)^rP_n^{(m)}}\rangle=\dfrac{1}{d^{(m)}_n}\begin{vmatrix}
\langle{u_1^{(m)},P_{n-m}}\rangle&\cdots&\langle{u_m^{(m)},P_{n-m}}\rangle&\langle{u_{p},(x-a)^{r}P_{n-m}}\rangle\\[5pt]
\vdots&&\vdots&\vdots\\
\langle{u_1^{(m)},P_{n-1}}\rangle&\cdots&\langle{u_m^{(m)},P_{n-1}}\rangle&\langle{u_{p},(x-a)^{r}P_{n-1}}\rangle\\[5pt]
\langle{u_1^{(m)},P_{n}}\rangle&\cdots&\langle{u_m^{(m)},P_{n}}\rangle&\langle{u_{p},(x-a)^{r}P_{n}}\rangle\\
\end{vmatrix}=0
\end{equation}
also because the last column of the determinant has all the entries equal to zero. 

Moreover, substituting $n$ by $ds+k-1,$ $s\in\mathbb{N}$ in \eqref{(b)}-\eqref{(c)} and expanding the determinants by their last columns we see
$$\langle{u_k^{(m)},(x-a)^rP_{dr+k-1}^{(m)}}\rangle=\begin{cases}
(-1)^{m}\dfrac{d^{(m)}_{n+1}}{d^{(m)}_n}\langle{u_{d-m+k},(x-a)^{r-1}P_{dr+k-1-m}}\rangle,\quad  k\leq m, \\[0.7cm]
(-1)^{m}\dfrac{d^{(m)}_{n+1}}{d^{(m)}_n}\langle{u_{p},(x-a)^{r}P_{dr+p-1}}\rangle, \quad  k=m+p,\ \ p=1,\ldots , d-m\,.
\end{cases}$$ 
Thus, 
\begin{equation}
\label{(d)}
\langle{u_k^{(m)},(x-a)^rP_{n-1}^{(m)}}\rangle\neq 0\,.
\end{equation}

From \eqref{(b)}-\eqref{(d)} we conclude that $(u^{(m)}_1,\ldots,u^{(m)}_{d})$ is a vector of orthogonality and $\{P^{(m)}_n\}$ is the corresponding monic $d$-OPS as we needed to prove.
$\hfill\square$

\subsection{Proof of Theorem \ref{principal}}

Now, we will see that the proof of Theorem \ref{principal} is an easy consequence of Section \ref{section2}. 

For $k=0,1,\ldots, d-1$, taking $r=d-k-1$ and $q=k+1$ in \eqref{nota2},
$$
(x-a){\bf P}^{(d-k-1)}(x)=N^{(d-k-1,k+1)}{\bf P}^{(d)}(x)\,.
$$
In particular,
\begin{equation}
\label{primero}
(x-a){\bf P}(x)=U{\bf P}^{(d)}(x)\,,
\end{equation}
where 
$U=N^{(0,d)}$
is a bi-diagonal upper triangular matrix whose entries in the diagonal (denoted by $s_0,\,s_1,\ldots$) are different from zero. That is, \eqref{primero} can be rewritten as
$$
(x-a)P_n(x)=P_{n+1}^{(d)}(x)+s_{n} P_n^{(d)}(x)\,,\quad n=0,1,\ldots.
$$
As a consequence  we have
$$
P^{(d)}_{n}(a)\neq 0\,,\quad n\in \mathbb{N}\,,
$$
because
\begin{equation}
\label{sexto}
P_{n+1}^{(d)}(a)+s_{n} P_n^{(d)}(a)=0
\end{equation}
and, if it were $P^{(d)}_{m}(a)=0$ for some $m\in \mathbb{N}$, then $P^{(d)}_{n}(a)=0$ for all $n\in \mathbb{N}$, which is incompatible with the recurrence relation \eqref{recurrence}.Thus,
$$
U=
\left(
\begin{array}{cccc}
-\frac{P_1^{(d)}(a)}{P_0^{(d)}(a)}& 1 \\
& -\frac{P_2^{(d)}(a)}{P_1^{(d)}(a)}& 1 \\
 & & \ddots  & \ddots  \\
\end{array}
\right). 
$$

With our hypotheses, we have the factorization \eqref{nuevo23} for $L^{(q,r)}$. Then, if $r+q=d$, we have
\begin{equation}
\label{(A)}
L^{(r,d-r)}=L^{(d)} L^{(d-1)}\cdots L^{(r+1)}.
\end{equation}

On the other hand, from \eqref{nota2} (also for $q=d-r$),
\begin{equation}
\label{(B)}
(x-a){\bf P}^{(r)}(x)=N^{(r,d-r)}{\bf P}^{(d)}(x).
\end{equation}
If $r=0$ then \eqref{(B)} becomes \eqref{primero}. If $r>0$, substituting $r$ by $r-1$ in \eqref{(B)} and multiplying by $L^{(r-1,1)}$ we obtain
$$
(x-a)L^{(r-1,1)}{\bf P}^{(r-1)}(x)=L^{(r-1,1)}N^{(r-1,d-r+1)}{\bf P}^{(d)}(x).
$$
From this and \eqref{nota1}, using the notation introduced in \eqref{otronuevo22},
\begin{equation}
\label{(C)}
(x-a){\bf P}^{(r)}(x)=L^{(r)}N^{(r-1,d-r+1)}{\bf P}^{(d)}(x).
\end{equation}
Comparing \eqref{(B)} and \eqref{(C)},
$$
N^{(r,d-r)}=L^{(r)}N^{(r-1,d-r+1)},
$$
which is verified for $r=1,2,\ldots,d$. 

Iterating the procedure,
\begin{eqnarray}
\label{(D)}
N^{(r,d-r)}=L^{(r)}L^{(r-1)}N^{(r-2,d-r+2)}=\cdots&= & L^{(r)}L^{(r-1)}\ldots  L^{(1)}N^{(0,d)}\nonumber\\
&= & L^{(r)}L^{(r)}\ldots  L^{(1)}U\,.
\end{eqnarray}
Finally, using \eqref{(A)} and \eqref{(D)} in \eqref{J} we arrive to \eqref{otrooctavo} and the result is proved.
$
\hfill \square
$

\end{document}